\documentclass[a4paper,12pt]{article}
\usepackage[a4paper,top=4.0cm,bottom=4.0cm,left=2.0cm,right=2.0cm]{geometry}

\usepackage{natbib}
\bibpunct{(}{)}{;}{d}{,}{,}

\usepackage{amsthm,amsmath,amsfonts,mathrsfs}
\usepackage{epsfig,graphicx}

\title{A localization technique for ensemble Kalman filters\thanks{Universit\"at Potsdam, 
Institut f\"ur Mathematik, Am Neuen Palais 10, D-14469 Potsdam, Germany}}
\author{
Kay Bergemann
\and
Sebastian Reich}

\begin{document}

\maketitle

\begin{abstract} Ensemble Kalman filter techniques are widely used to assimilate observations into
dynamical models. The phase space dimension is typically much larger than the number of ensemble 
members which leads to inaccurate results in the computed covariance matrices.  These inaccuracies
can lead, among other things, to spurious long range correlations which can be eliminated by 
Schur-product-based localization techniques. In this paper, we propose a new technique for 
implementing such localization techniques within the class of ensemble transform/square root 
Kalman filters. Our approach relies on a continuous embedding of the Kalman filter update for the 
ensemble members, i.e., we state an ordinary differential equation (ODE) whose solutions, over a unit time
interval, are equivalent to the Kalman filter update. The ODE formulation forms a gradient system with
the observations as a cost functional. Besides localization, the new ODE ensemble formulation should also
find useful applications in the context of nonlinear observation operators and observations arriving 
continuously in time.
\end{abstract}

\noindent
{\bf Keywords.}
Data assimilation, ensemble Kalman filter, localization, continuous Kalman filter


\section{Introduction}

We consider ordinary differential equations
\begin{equation} \label{ode}
\dot{\bf x} = f({\bf x},t)
\end{equation}
with state variable ${\bf x} \in \mathbb{R}^n$. Initial conditions at time $t_0$ are not 
precisely known and we assume instead that
\begin{equation} \label{IC}
{\bf x}(t_0) \sim {\rm N}({\bf x}_0,{\bf B}),
\end{equation}
where ${\rm N}({\bf x}_0,{\bf B})$ denotes an $n$-dimensional Gaussian distribution with mean
${\bf x}_0 \in \mathbb{R}^n$ and covariance matrix ${\bf B} \in \mathbb{R}^{n \times n}$. We also
assume that we obtain measurements ${\bf y}(t_i) \in \mathbb{R}^{k}$ at discrete times
$t_j\ge t_0$, $j=0,1,\ldots,M$, subject to measurement errors, which are also Gaussian 
distributed with zero mean and covariance matrix ${\bf R} \in \mathbb{R}^{k\times k}$, i.e.,
\begin{equation} \label{measurement}
{\bf y}(t_j) - {\bf H}{\bf x}(t_j) \sim {\rm N}({\bf 0},{\bf R}).
\end{equation}
Here ${\bf H} \in \mathbb{R}^{k\times n}$ is the (linear) measurement operator.

Data assimilation is the task to combine the model (\ref{ode}) (here assumed to be perfect), 
the knowledge about the initial conditions (\ref{IC}) and available measurements 
(\ref{measurement}) in a prediction of the probability distribution of the solution
at any time $t>t_0$. We refer to \cite{sr:Lewis} for a detailed introduction and 
available approaches to data assimilation. In this paper, we
focus on the ensemble Kalman filter (EnKF) method, originally proposed by Evensen (see 
\cite{sr:evensen} for a recent account) and, in particular, on ensemble transform 
\citep{sr:bishop01}, ensemble adjustment \citep{sr:anderson01}, and 
ensemble square root filters \citep{sr:tippett03} and their sequential implementation
\citep{sr:whitaker02,sr:anderson03}.

The EnKF relies on the simultaneous propagation of $m$ independent solutions
${\bf x}_i(t)$, $i=1,\ldots,m$, from which we can extract an empirical mean
\begin{equation}
\overline{\bf x}(t) = \frac{1}{m} \sum_{i=1}^m {\bf x}_i(t)
\end{equation}
and an empirical covariance matrix
\begin{equation}
{\bf P}(t) = \frac{1}{m-1} \sum_{i=1}^m \left({\bf x}_i(t)-\overline{\bf x}(t)\right)
 \left({\bf x}_i(t)-\overline{\bf x}(t)\right)^T.
\end{equation}
In typical applications from meteorology, the ensemble size $m$ is much smaller than the
dimension $n$ of the model phase space and, more importantly, also much smaller than
the number of positive Lyapunov exponents.
Hence ${\bf P}(t)$ is highly rank deficient which can 
lead to unreliable predictions. Ensemble localization has been introduced by 
\cite{sr:houtekamer01} and \cite{sr:hamill01} to overcome this problem. However, only 
two techniques are currently available to implement Schur-product-based localization within 
the framework of ensemble transform/square root Kalman filters. The first option is provided by a 
sequential processing of observations \citep{sr:whitaker02,sr:anderson03}, 
while the deterministic ensemble Kalman filter (DEnKF) of \cite{sr:sakov08} 
is a second, more recent, option. 
The DEnKF  results in an approximate implementation of ensemble transform/square root Kalman filters. 
We also mention box/local analysis methods \citep{sr:evensen03,sr:ott04,sr:hunt07}, 
which assimilate data locally  in physical space and which therefore 
possess a ``built in'' localization. 

In this paper, we demonstrate that techniques proposed by \cite{sr:bgr09} 
for the filter analysis step can be further generalized to an ordinary differential 
equation (ODE) formulation in terms of the ensemble members ${\bf x}_i$, $i=1,\ldots,m$. 
This formulation is subsequently used to derive an easy to implement localized ensemble 
Kalman filter, which can process observations simultaneously and can be
extended to nonlinear observation operators.


\section{Background material}

We summarize a number of key results and techniques regarding ensemble Kalman filters.
We refer to \cite{sr:evensen} for an introduction and in-depth discussion of such
filters. 

\subsection{Kalman analysis step}

Let $n$ denote the dimension of the phase space of the problem. 
We consider an ensemble of $m$ members ${\bf x}_i(t) \in \mathbb{R}^n$ which we collect
in a matrix ${\bf X}(t) \in \mathbb{R}^{n\times m}$. In terms of ${\bf X}$, 
the ensemble mean is given by
\begin{equation}
\overline{\bf x}(t) = \frac{1}{m} {\bf X}(t) {\bf e}
\in \mathbb{R}^n
\end{equation}
and we introduce the ensemble deviation matrix 
\begin{equation}
{\bf X}'(t) = {\bf X}(t) - \overline{\bf x}(t) {\bf e}^T \in \mathbb{R}^{n\times m},
\end{equation}
where ${\bf e} = (1,\ldots,1)^T \in \mathbb{R}^m$.

We now describe the basic Kalman analysis step. Let $\overline{\bf x}_f$
and ${\bf X}'_f$ denote the forecast mean and deviation matrix, respectively.
The ensemble mean is updated according to
\begin{equation} \label{Kalmanmean}
\overline{\bf x}_a = \overline{\bf x}_f - {\bf K} \left( {\bf H}
\overline{\bf x}_f - {\bf y}\right),
\end{equation}
where 
\begin{equation} \label{Kalmangain}
{\bf K} = {\bf P}_f {\bf H}^T\left( {\bf H} {\bf P}_f {\bf H}^T + {\bf R}
\right)^{-1}
\end{equation}
is the Kalman gain matrix with empirical covariance matrix
\begin{equation}
{\bf P}_f = \frac{1}{m-1} {\bf X}'({\bf X}')^T,
\end{equation}
and ${\bf R} \in \mathbb{R}^{k\times k}$ is the measurement error covariance matrix.

While the update of the mean is common to most ensemble Kalman filters, 
the update of the ensemble deviation matrix ${\bf X}'_f$ can be implemented in
several ways. In this paper, we focus on ensemble update techniques that
employ either a transformation of the form
\begin{equation} \label{transform1}
{\bf X}'_a = {\bf A} {\bf X}'_f
\end{equation}
with an appropriate matrix ${\bf A} \in \mathbb{R}^{n\times n}$ \citep{sr:anderson01} or
a transformation
\begin{equation} \label{transform2}
{\bf X}'_a = {\bf X}'_f {\bf T}
\end{equation}
with an appropriate ${\bf T} \in \mathbb{R}^{m\times m}$ 
\citep{sr:bishop01,sr:whitaker02,sr:tippett03,sr:evensen04}. 
The matrices ${\bf A}$ and 
${\bf T}$ are chosen such that the resulting ensemble deviation matrix 
${\bf X}'_a$ satisfies
\begin{equation} \label{Kalmanupdate2}
{\bf P}_a = \frac{1}{m-1}{\bf X}'_a ({\bf X}'_a)^T = \left({\bf I} - {\bf K}
{\bf H}\right) {\bf P}_f.
\end{equation}
It has been shown by \cite{sr:tippett03} that
both formulations (\ref{transform1}) and (\ref{transform2}) can be made equivalent not only in terms of
(\ref{Kalmanupdate2}) but also in terms of the resulting ensemble deviation matrix ${\bf X}_a'$. 
Since $n\gg m$ in most applications, formulation (\ref{transform2}) is generally preferred except
when working in a sequential framework \citep{sr:anderson03}. 

Note that the transformation matrix ${\bf T}$ should also satisfy ${\bf T}{\bf e} = {\bf e}$ 
to guarantee
${\bf X}'_a {\bf e} = {\bf 0}$ \citep{sr:wang04,sr:nichols08,sr:sakov08b}. 
Otherwise, the update of the ensemble
deviation matrix would affect the update of the ensemble mean.

Several methods have been proposed recently (including \cite{sr:sakov08} and
\cite{sr:bgr09}) that satisfy (\ref{Kalmanupdate2}) only approximately. More specifically, 
\cite{sr:sakov08} suggest to use (\ref{transform1}) with
\begin{equation} \label{denkf}
{\bf A}  = {\bf I} - \frac{1}{2} {\bf K} {\bf H},
\end{equation}
while \cite{sr:bgr09} use numerical approximations to the underlying ODE formulation
\begin{equation} \label{continuousKalman1}
\frac{{\rm d}}{{\rm d} s}{\bf Y} = -\frac{1}{2m-2}
{\bf Y} {\bf Y}^T {\bf H}^T {\bf R}^{-1} {\bf H} {\bf Y}
\end{equation}
in a fictitious time $s \in [0,1]$
See, for example, \cite{sr:simon} for a derivation of (\ref{continuousKalman1}).
The initial condition is ${\bf Y}(0) = {\bf X}'_f$ and the updated ensemble deviation matrix,
which satisfies (\ref{Kalmanupdate2}) exactly, is provided by the solution
at time $s=1$, i.e.
\begin{equation}
{\bf X}'_a = {\bf Y}(1).
\end{equation}
A typical numerical implementation of (\ref{continuousKalman1}) uses two 
or four time-steps with the forward Euler method \citep{sr:bgr09}. The resulting
transformation of the forecast into the analyzed ensemble
deviation matrix is of the form (\ref{transform1}) with ${\bf A}$ defined through
the time-stepping method.

Note that the Kalman gain matrix (\ref{Kalmangain}) is equivalent to
\begin{equation} \label{Kalmangain2}
{\bf K} = {\bf P}_a {\bf H}^T {\bf R}^{-1},
\end{equation}
which is advantageous in connection with (\ref{continuousKalman1}) since
only the inversion of the measurement error covariance matrix 
${\bf R} \in \mathbb{R}^{k\times k}$
is now required to implement a complete Kalman analysis step. Algorithmically,
one would first update the ensemble deviation matrix using 
(\ref{continuousKalman1}), then form the analysed ensemble covariance matrix
${\bf P}_a = {\bf X}_a'[{\bf X}_a']^T/(m-1)$ as well as the Kalman gain
matrix (\ref{Kalmangain2}), and finally update the ensemble mean using
(\ref{Kalmanmean}).

All methods discussed so far have in common that the Kalman update increments for the ensemble mean 
and ensemble deviation matrix lie in a $m-1$ dimensional
subspace, denoted by $\mathbb{S}_f \subset \mathbb{R}^n$. This space is defined by 
the range/image of the forecast ensemble deviation matrix 
${\bf X}_f'$. \cite{sr:bgr09} introduced a continuous matrix factorization algorithm for
the ensemble ${\bf X}(t)$, which automatically produces orthogonal 
vectors that span $\mathbb{S}_f$. 

It is common practice  
to apply variance inflation 
\citep{sr:andand99} to ${\bf X}'(t_j)$ before the forecasted ensemble 
is updated under the Kalman filter analysis step, i.e., the Kalman analysis
step uses
\begin{equation} \label{inflation}
{\bf X}_f := \overline{\bf x}(t_j){\bf e}^T + \delta\, {\bf X}'(t_j),
\end{equation}
where $\delta \ge 1$ is an inflation factor, instead of
${\bf X}_f = {\bf X}(t_j)$.

\subsection{Localization}

The idea of localization, as proposed by \cite{sr:houtekamer01}, 
is to replace the matrices ${\bf H} {\bf P}_f$ and
${\bf H}{\bf P}_f {\bf H}^T$ in the Kalman gain matrix (\ref{Kalmangain}) by
\begin{equation}
\widetilde{{\bf H} {\bf P}_f } = {\bf C}_{{\rm loc},1} \circ 
\left({\bf H} {\bf P}_f \right), \quad
\widetilde{{\bf H} {\bf P}_f {\bf H}^T} = {\bf C}_{{\rm loc},2} \circ \left({\bf H} 
{\bf P}_f {\bf H}^T\right),
\end{equation}
respectively, where ${\bf C}_{{\rm loc},1} \in \mathbb{R}^{n\times k}$ and 
${\bf C}_{{\rm loc},2} \in \mathbb{R}^{k\times k}$
are appropriate localization matrices based on filter functions suggested by 
\citep{sr:gaspari99} and ${\bf C} \circ {\bf Y}$ denotes the Schur
product of two matrices ${\bf C}$ and ${\bf Y}$ of identical dimension, i.e.,
\begin{equation}
\left({\bf C} \circ {\bf Y}\right)_{i,j} = ({\bf C})_{i,j} \,({\bf Y})_{i,j}
\end{equation}
for all indices $i,j$. We denote the resulting modified Kalman gain matrix
by ${\bf K}_{\rm loc,f}$, i.e., 
\begin{equation} \label{Kalmangainloc1}
{\bf K}_{\rm loc,f} = (\widetilde{{\bf H} {\bf P}_f})^T \left(
\widetilde{{\bf H} {\bf P}_f {\bf H}^T} + {\bf R}\right)^{-1}.
\end{equation}
Localization was also proposed by \cite{sr:hamill01} with the only difference that
localization is not applied to ${\bf H} {\bf P}_f {\bf H}^T$. 

Alternatively, one can localize the Kalman gain matrix formulation (\ref{Kalmangain2}) and use
\begin{equation} \label{Kalmangainloc2}
{\bf K}_{\rm loc,a} = (\widetilde{{\bf H} {\bf P}_a })^T  {\bf R}^{-1}
\end{equation}
instead of (\ref{Kalmangain}) in an ensemble Kalman filter.
Note that (\ref{Kalmangainloc1}) and (\ref{Kalmangainloc2}) are not equivalent in general 
and that formulation (\ref{Kalmangainloc2}) is easier to implement.

Based on these modified Kalman gain matrices, a Schur-product-based
localization is easy to apply to the update (\ref{Kalmanmean}) of the ensemble mean and
to ensemble deviation updates that use perturbed observations \citep{sr:burgers98}, which
is essentially the localization approach of \cite{sr:houtekamer01} and \cite{sr:hamill01}.
However, the popular class of ensemble transform/square root
filters, based on (\ref{transform2}), has not yet been amenable to
Schur-product-based localizations except when observations are treated sequentially, i.e., 
when $k=1$ in each transformation step \citep{sr:whitaker02}. 

It is feasible that localizations can be implemented for ensemble deviation updates of the form 
(\ref{transform1}) through an appropriate modification of the ensemble adjustment technique of 
\cite{sr:anderson01}. However, such a modification would lead to a computationally expensive 
implementation of Schur-product-based localizations. The recently proposed DEnKF filter of 
\cite{sr:sakov08}, on the other hand, leads to a computationally feasible implementation with
the localization directly applied to (\ref{denkf}), i.e., one uses
\begin{equation}
{\bf A} = {\bf I} - \frac{1}{2} {\bf K}_{\rm loc,f} {\bf H}
\end{equation}
in (\ref{transform1}).

We note that localization implies in general that the update increments for the ensemble mean
and the ensemble deviation matrix no longer lie in the subspace
$\mathbb{S}_f$ defined by the range/image of ${\bf X}_f'$. While this is a desirable
property on the one hand, it can lead to unbalanced fields in the 
analyzed ensemble ${\bf X}_a$ on the other hand. This has been investigated, 
for example, by \cite{sr:houtekamer05} and \cite{sr:kepert09}. 

We finally mention an alternative approach to localization. 
The box/local EnKF  filters of \cite{sr:evensen03,sr:ott04,sr:hunt07} 
assimilate data locally in physical space and possess a
``built in'' localization based on the spatial structure of the underlying partial 
differential equation model.


\section{Localization based on continuous ensemble updates}

We now describe an alternative for introducing localization, which is based on 
a generalization of the ODE formulation (\ref{continuousKalman1}) and which
leads to an ODE formulation directly in the ensemble members ${\bf x}_i$. 

We first note that the Kalman update (\ref{Kalmanmean}) for the ensemble mean can also be formulated
in terms of an ODE, i.e.,
\begin{equation} \label{continuousKalman2}
\frac{\rm d}{{\rm d} s} \overline{\bf x} = - 
\frac{1}{m-1} {\bf Y} {\bf Y}^T {\bf H}^T {\bf R}^{-1} \left(
{\bf H} \overline{\bf x} - {\bf y}\right) 
\end{equation}
with $\overline{\bf x}(0) = \overline{\bf x}_f$ and 
$\overline{\bf x}_a = \overline{\bf x}(1)$. See, for example, 
\cite{sr:simon} for a derivation of (\ref{continuousKalman2}).

To further reveal the underlying mathematical structure of
(\ref{continuousKalman1}) and (\ref{continuousKalman2}), 
we introduce the cost functional
\begin{equation} \label{cost1}
S({\bf x}) = \frac{1}{2} \left({\bf H}{\bf x} -{\bf y}\right)^T
{\bf R}^{-1} \left({\bf H}{\bf x} -{\bf y}\right)
\end{equation}
for each set of observations. Next we 
combine (\ref{continuousKalman1}) and (\ref{continuousKalman2})
to give rise to the differential equations
\begin{equation} \label{continuousKalman3}
\frac{\rm d}{{\rm d}s} {\bf x}_i = 
- \frac{1}{2} {\bf P} \left\{ \nabla_{{\bf x}_i} S({\bf x}_i) + 
\nabla_{\overline{\bf x}} S(\overline{\bf x})\right\}
\end{equation}
in the ensemble members ${\bf x}_i$, $i=1,\ldots,m$. The equations
are closed through the standard definition 
\begin{equation}
\overline{\bf x}(s) = \frac{1}{m} \sum_{i=1}^m {\bf x}_i(s)
\end{equation}
for the mean and covariance matrix
\begin{equation}
{\bf P}(s) = \frac{1}{m-1} \sum_{i=1}^m \left({\bf x}_i(s) - 
\overline{\bf x}(s)\right) \left({\bf x}_i(s) - 
\overline{\bf x}(s)\right)^T.
\end{equation}
Since the covariance matrix ${\bf P}$ is symmetric, a straightforward calculation 
reveals that
\begin{equation}
\frac{\rm d}{{\rm d}s} \left\{
S(\overline{\bf x}) + \frac{1}{m} \sum_{i=1}^m S({\bf x}_i) 
\right\} \le 0
\end{equation}
along solutions of (\ref{continuousKalman3}). More precisely, (\ref{continuousKalman3})
is equivalent to the gradient system 
\begin{equation} \label{gradient}
\frac{\rm d}{{\rm d}s} {\bf x}_i =  -{\bf P} \nabla_{{\bf x}_i} {\cal V}({\bf X}),
\end{equation}
in the ensemble matrix ${\bf X}(s)$ with potential
\begin{equation} \label{potential}
{\cal V}({\bf X}) = \frac{m}{2} \left\{ S(\overline{\bf x}) + \frac{1}{m} \sum_{i=1}^m S({\bf x}_i)
\right\} .
\end{equation}
The actual decay of the potential ${\cal V}$ along solutions of (\ref{gradient}) 
depends crucially on the covariance matrix ${\bf P}$.

We note that Schur-product-based localizations can easily be applied
to (\ref{continuousKalman3}) to obtain, e.g., 
\begin{equation} \label{continuousKalman4a}
\frac{\rm d}{{\rm d}s} {\bf x}_i = 
- \frac{1}{2} \widetilde{\bf P} \left\{ \nabla_{{\bf x}_i} S({\bf x}_i) + 
\nabla_{\overline{\bf x}} S(\overline{\bf x})\right\}, \qquad \widetilde{\bf P} = {\bf C}_{\rm loc} \circ {\bf P},
\end{equation}
or
\begin{equation} \label{continuousKalman4b}
\frac{\rm d}{{\rm d}s} {\bf x}_i = 
- \frac{1}{2} (\widetilde{\bf H P})^T {\bf R}^{-1} \left\{  {\bf H}{\bf x}_i + {\bf H}\overline{\bf x}
- 2{\bf y}  \right\}, \qquad \widetilde{\bf H P} = {\bf C}_{{\rm loc},1} \circ ({\bf H P}),
\end{equation}
in case of linear observation operators. These modified ensemble update formulations 
are easy to implement numerically. See Section \ref{sec_implementation}
for details.


\section{Numerical implementation aspects} \label{sec_implementation}

The various ODE formulations for the ensemble members ${\bf x}_i$, $i=1,\ldots,m$, 
need to be solved over a unit time interval with initial conditions provided by the forecast 
values ${\bf x}_{i,f}$ of the ensemble members.  We apply the forward Euler
method with step-size $\Delta s = 1/4$ (four time-steps) for our experiments. We found that 
$\Delta s = 1$ (single time-step) and $\Delta s = 1/3$ (three time-steps) lead to unstable 
simulations, while $\Delta s  = 1/2$  (two time-steps) leads to occasional instabilities for
larger values of the ensemble inflation factor $\delta$ in (\ref{inflation}). 
On the other hand, increasing the number of time-steps beyond four did not change the results 
significantly. We also expect that four time-steps will generally be sufficient in practical 
applications unless observations strongly contradict their forecast values and large gradient
values are generated in (\ref{continuousKalman3}). The same consideration can apply to simulations with
large inflation factors $\delta$ in (\ref{inflation}). As a safe guard, one can monitor the decay
of the potential (\ref{potential}) along numerically generated solutions and adjust the step-size
$\Delta s$ if necessary.

Note that the continuous formulations do not require matrix inversions/factorizations except
for the computation of ${\bf R}^{-1}$. The computational cost of localization can be reduced even 
further by using the following approximation.
The matrix $\widetilde{\bf HP}$ in (\ref{continuousKalman4b}) varies along
solutions and an approximative formulation is obtained by replacing
$\widetilde{{\bf H}{\bf P}}(s)$ by its value at $s =0$ for all $s>0$. This leads
to a linear ODE in the ensemble members ${\bf x}_i$ with constant coefficient matrix, i.e.,
\begin{equation} \label{continuousKalman4c}
\frac{{\rm d}}{{\rm d} s}{\bf x}_i = - \frac{1}{2} (\widetilde{{\bf H}{\bf P}}(0))^T 
{\bf R}^{-1}  \left\{  {\bf H}{\bf x}_i + {\bf H}\overline{\bf x}
- 2{\bf y}  \right\} .
\end{equation}
Note that $\widetilde{\bf H P}$ and $\widetilde{{\bf HPH}^T}$ are sparse matrices 
for compactly supported filter functions \citep{sr:gaspari99}. Numerical implementations 
of (\ref{continuousKalman4c}) should first update the increments 
${\bf z}_i = {\bf H}{\bf x}_i-{\bf y}$ 
with Euler's method, i.e.,
\begin{equation}
{\bf z}_i^{l+1} = {\bf z}_i^l - \frac{\Delta s}{2} \widetilde{ {\bf H P H}^T}(0) 
{\bf R}^{-1} \left\{ {\bf z}_i^l + \frac{1}{m} \sum_{j=1}^m {\bf z}_j^l\right\}, \qquad 
l = 0,\ldots,L,
\end{equation}
$L = 1/\Delta s$ the number of integration steps,
and then use the accumulated increments 
\begin{equation}
{\bf z}_i = \sum_{l=0}^{L-1} {\bf z}_i^l
\end{equation} 
to compute the final
update of the ensemble members ${\bf x}_i$, $i=1,\ldots,m$. Overall matrix-vector-products will 
induce a computational complexity of ${\cal O}(k\,m)$ 
in the ensemble size $m$ and the number of observations $k$ 
independent of the system size $n$. The same order of complexity applies to the serial 
algorithm of \cite{sr:hamill01} with the important difference that (\ref{continuousKalman4c}) can
be implemented as a simultaneous update over all observations. 

\cite{sr:bgr09} proposed a re-orthogonalization technique for the ensemble deviation matrix
${\bf X}'$. It should be noted that the re-orthogonalization is not
uniquely defined. We implemented several variants of re-orthogonalization in combination with
localization but did not find any significant improvements in the results.


\begin{figure}
\begin{center}
\includegraphics[width=0.5\textwidth]{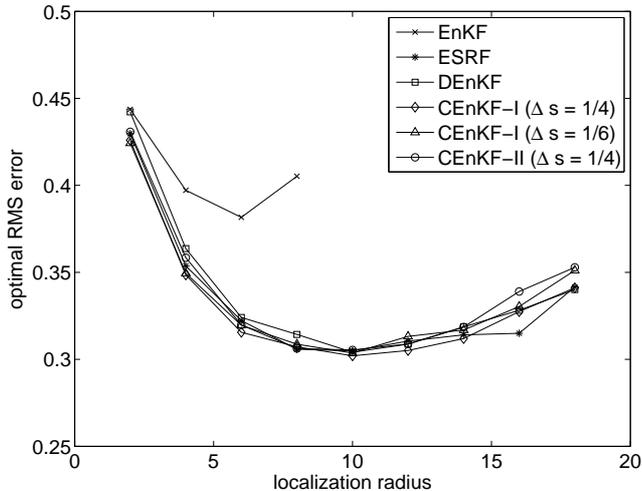}
\end{center}
\caption{The best RMS error for the Lorenz-96 model with an ensemble size
of $m = 10$ and $k=20$ observations taken in intervals of $\Delta t_{\rm obs} = 0.05$ over a total
of 5000 assimilation cycles.}
\label{fig_lorenz1}
\end{figure}


\section{Numerical experiments} \label{sec_numerics}

We now report results from two test problems 
and implementations of  (\ref{continuousKalman4b}) and
(\ref{continuousKalman4c}) with localization. The results are compared to those from
standard ensemble Kalman filter techniques.

\subsection{Lorenz-96 model}

The standard implementation of the Lorenz-96 model \citep{sr:lorenz96,sr:lorenz98} has state 
vector ${\bf x} = (x_1,\ldots,x_{n})^T \in \mathbb{R}^{n}$, $n=40$, 
and its time evolution is given by the differential equations
\begin{equation}
\dot{x}_j = (x_{j+1}-x_{j-2})x_{j-1} - x_j + 8
\end{equation}
for $j=1,\ldots,n$. To close the equations, we define 
$x_{-1} = x_{39}$, $x_0 = x_{40}$,
and $x_{41} = x_1$. 

The attractor of this standard implementation has a fractal dimension of about 27 and 13 positive
Lyapunov exponents. Localization will be necessary for ensembles with $m\le 13$ 
ensemble members. We use an ensemble size of 
$m=10$ in our experiments. We observe every second grid point, i.e., $k=20$, and the measurement error
covariance is ${\bf R} = {\bf I}_{20}$. Measurement are taken in time intervals of 
$\Delta t_{\rm obs} = 0.05$. After a short spin-up period, a total of $J = 5000$ analysis steps are performed
in each experiment. The "true" trajectory ${\bf x}_{\rm truth}(t_n)$ 
is generated by integrating the Lorenz-96 model with the implicit midpoint rule and step-size
$\Delta t = 0.005$, i.e., we assume that there is no model error. 
The observations are obtained according to
\begin{equation}
{\bf y}(t_{\rm obs}) = {\bf H}{\bf x}_{\rm truth}(t_{\rm obs}) 
+ {\bf r}(t_{\rm obs}),
\end{equation}
where ${\bf r}(t_{\rm obs})$ are i.i.d.~Gaussian random numbers with
mean zero and covariance matrix ${\bf R}$.

We implement localization combined with standard ensemble inflation 
for the following five different ensemble Kalman filters:
(i) EnKF with perturbed observations \citep{sr:burgers98,sr:houtekamer05},
(ii) ensemble square root filter (ESRF) 
with sequential treatment of observations  \citep{sr:whitaker02},
(iii) DEnKF \citep{sr:sakov08},
(iv) formulation (\ref{continuousKalman4b}), denoted CEnKF-I,
(v) formulation (\ref{continuousKalman4c}), denoted CEnKF-II.
We implement CEnKF-I with $\Delta s = 1/4$ and $\Delta s = 1/6$, respectively, 
to demonstrate the impact of the discretization parameter
on the results. 

For simplicity,
localization is performed by multiplying each element of the matrices ${\bf H} {\bf P}$ and
${\bf H} {\bf P} {\bf H}^T$, respectively, by a distance dependent factor
$\rho_{i,i'}$. This factor is defined by the compactly supported localization function
(4.10) from \cite{sr:gaspari99}, distance 
$r_{i,i'} = \min\{|i-i'|,n-|i-i'|\}$, where $i$ and $i'$ denote the
indices of the associated observation/grid points $x_{i}$ and $x_{i'}$, respectively, and
a fixed localization radius $r_0$. The localization radius is varied between $r_0 = 2$ and
$r_0 = 30$. The inflation factor $\delta$ in (\ref{inflation}) is taken from the range
$[\sqrt{1.02},\sqrt{1.16}]$.

In Figure \ref{fig_lorenz1}, we display the RMS error 
\begin{equation} \label{rmse}
{\rm rmse} = \sqrt{ \frac{1}{n\,J} \sum_{j=1}^J \|  \overline{\bf x}(j\cdot \Delta t_{\rm obs}) - 
 {\bf x}_{\rm truth}(j\cdot \Delta t_{\rm obs})\|^2 }
\end{equation}
for an optimally chosen inflation factor $\delta$ as a function
of the localization radius $r_0$. Results are displayed for those 
localization radii $r_0$ only, which lead to at least one simulation with 
a RMS error of less than one.

We conclude that EnKF yields the lowest filter skills while all other
methods show an almost identical performance.


\begin{table} 
\begin{center}
 {\scriptsize
\begin{tabular}{|c||c|c|c|c|c|c|c|}
\hline
$\delta \backslash r_0$ & {\bf 5}& {\bf 10}& {\bf 15}& {\bf 20}& {\bf 25}& {\bf 30}& {\bf 35}\\ \hline\hline
{\bf  1.02} &  0.59 &  0.62 &  0.72 &  0.96 &  1.41 & Inf  & Inf \\ \hline
{\bf  1.06} &  0.80 &  0.67 &  0.62 &  0.65 &  0.79 & 0.99 & 1.66 \\ \hline
{\bf  1.10} &  1.09 &  0.85 &  0.72 &  0.69 &  0.75 & 1.04 & 1.13 \\ \hline
{\bf  1.14} & 1.38  &  1.01 &  0.84 &  0.77 &  0.78 & 0.88 & 1.40 \\ \hline
{\bf  1.18} & Inf   &  1.18 &  0.96 &  0.84 &  0.82 & 0.86 &  1.35 \\ \hline
\end{tabular}  
} 
\end{center}

\caption{Mean RMS error for localized CEnKF-I over 4000 time steps as a functions
of the localization radius $r_0$ and the inflation factor $\delta$. For clarity, the value Inf 
is assigned if the RMS error exceeds the value 2.0 (no filter skill).}
\label{table1}

\begin{center}
 {\scriptsize
\begin{tabular}{|c||c|c|c|c|c|c|c|}
\hline
$\delta \backslash r_0$ & {\bf 5}& {\bf 10}& {\bf 15}& {\bf 20}& {\bf 25}& {\bf 30}& {\bf 35} \\ \hline\hline
{\bf  1.02} &  0.60 &  0.65 &  0.77 &  1.01 &  1.31 & 1.80 & Inf  \\ \hline
{\bf  1.06} &  0.80 &  0.66 &  0.63 &  0.68 &  0.86 &  1.13 & 1.51  \\ \hline
{\bf  1.10} &  1.11 &  0.84 &  0.72 &  0.69 &  0.74 &  0.93 &  1.29 \\ \hline
{\bf  1.14} &  1.42 &  1.04 &  0.83 &  0.76 &  0.78 &  0.86 &  1.14 \\ \hline
{\bf  1.18} & Inf  &  1.21  &  0.96 &  0.85 &  0.82 & 0.86 &  1.16  \\ \hline
\end{tabular}  
} 
\end{center}

\caption{Mean RMS error for localized CEnKF-II
over 4000 time steps as a functions
of the localization radius $r_0$ and the inflation factor $\delta$. 
}
\label{table2}
\end{table}

\begin{table} 
\begin{center}
{\scriptsize
\begin{tabular}{|c||c|c|c|c|c|c|c|}
\hline
$\delta \backslash r_0$ & {\bf 5}& {\bf 10}& {\bf 15}& {\bf 20}& {\bf 25}& {\bf 30}& {\bf 35} \\ \hline\hline
{\bf  1.02} &  0.59 &  0.62 &  0.75 &  0.94 &  1.06 & Inf & Inf \\ \hline
{\bf  1.06} &  0.82 &  0.69 &  0.64 &  0.68 &  0.73 & 0.98  & Inf \\ \hline
{\bf  1.10} &  1.16 &  0.89 &  0.75 &  0.70 &  0.72 & 0.85 &  1.22 \\ \hline
{\bf  1.14} &  1.50 &  1.11 &  0.89 &  0.80 &  0.77 & 0.85 &  0.99 \\ \hline
{\bf  1.18} & Inf  &  1.33 &  1.05 &  0.91 &  0.87 & 0.87 &  1.01 \\ \hline
\end{tabular} 
}
\end{center}
\caption{Mean RMS error for localized DEnKF over 4000 time steps as a functions
of the localization radius $r_0$ and the inflation factor $\delta$. 
}
\label{table3}
\end{table}

\begin{figure}
\begin{center}
\includegraphics[width=0.5\textwidth]{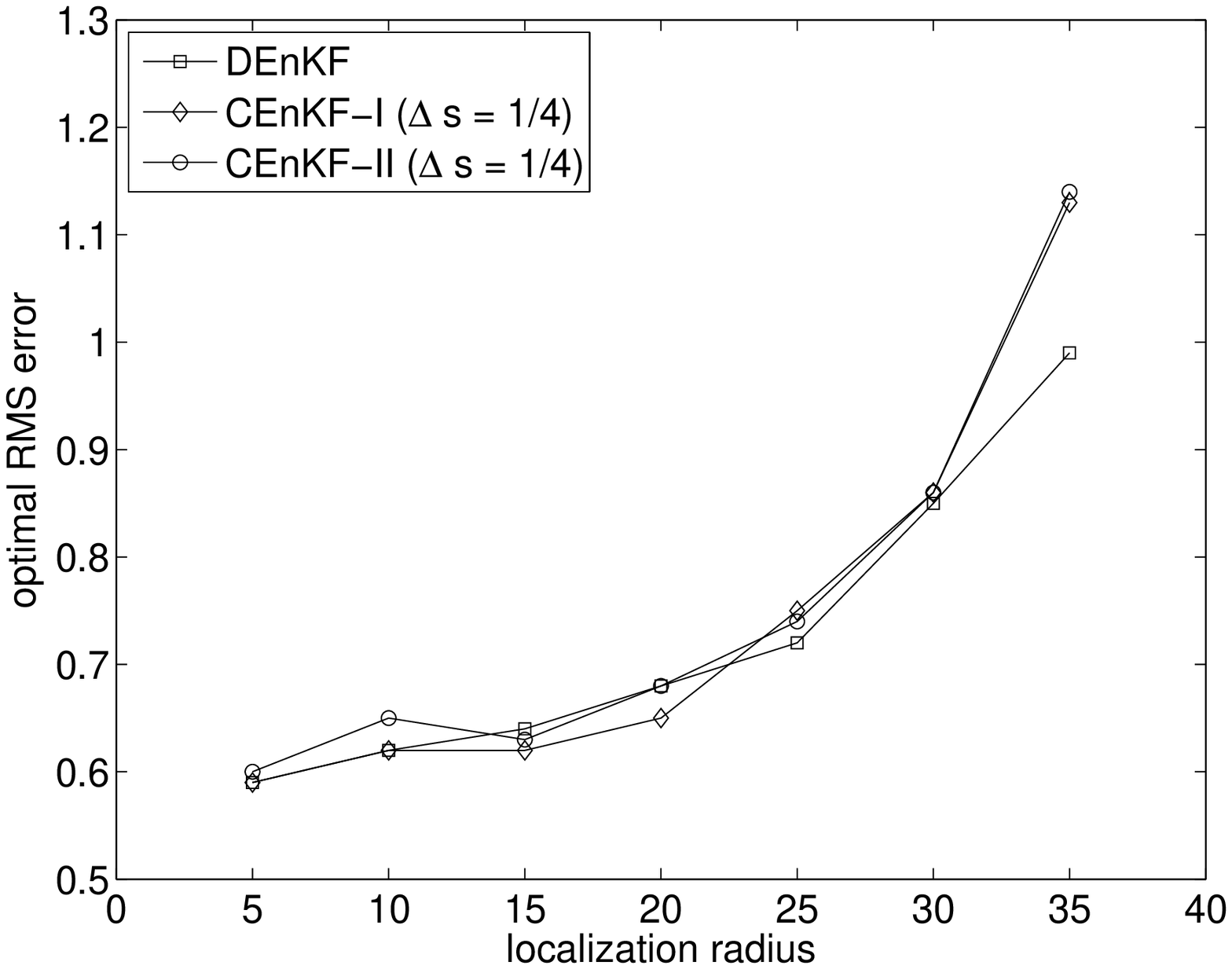}
\end{center}
\caption{The best RMS error for the QG-model of \cite{sr:sakov08} with an ensemble size
of $m = 25$ and $k=300$ observations taken in intervals of $\Delta t_{\rm obs} = 4.0$ over a total
of 1000 assimilation cycles.}
\label{fig_qg1}
\end{figure}


\subsection{A quasi-geostrophic (QG) model}

We use the QG model of \cite{sr:sakov08}. The QG model is a numerical approximation of the
following 1.5-layer reduced gravity quasi-geostrophic model with double-gyre wind forcing and biharmonic
friction:
\begin{equation} \label{QG}
q_t = -\psi_x - \varepsilon J(\psi,q) - A \Delta ^3 \psi + 2\pi \sin (2\pi y),
\end{equation}
where $q = \Delta \psi - F\psi$, $J(\psi,q) = \psi_x q_y - \psi_y q_x$, $\Delta = \partial^2 /
\partial x^2 + \partial^2 /\partial y^2$. The coefficients are given by $F = 1600$, $\varepsilon = $
10$^{-5}$, $A = 2 \times $10$^{-12}$. The model domain is $(x,y) \in [0,1]\times [0,1]$ 
with zero Dirichlet boundary conditions. The model is discretized
over this domain using a $129 \times 129$ grid. For more details 
see \cite{sr:sakov08}.

We implement the deterministic ensemble Kalman filter 
(DEnKF) and our ensemble Kalman filters based on (\ref{continuousKalman4b}). All experiments
use $m=25$ ensemble members. The dimension of the phase space is $n= 16129$. The dimension of
the attractor and the number of positive Lyapunov exponents are currently not known.

In line with \cite{sr:sakov08},
localization is performed by multiplying each element of the matrices ${\bf H} {\bf P}$ and
${\bf H} {\bf P} {\bf H}^T$, respectively, by a factor $\rho_{ij,i'j'} = 
\exp(-0.5 r_{ij,i'j'}^2/r_0^2)$. Here we use the distance 
$r_{ij,i'j'} = \sqrt{|i-i'|^2 + |j-j'|^2}$, 
where $(i,j)$ and $(i',j')$ denote the
indices of the associated observation/grid points $x_{ij}$ and $x_{i'j'}$, respectively, and
$r_0$ is a fixed localization radius. 

We test the performance of the filters for different values of the ensemble inflation factor $\delta$
in (\ref{inflation}) and the localization radius $r_0$. For each pair $(\delta,r_0)$ of simulation parameters, 
we run a single simulation over 4000 time steps with step-size $\Delta t = 1.25$ and 
perform a total of 1000 assimilation cycles using 300 observations of $\psi$ with observation variance
of 4.0 as described in \cite{sr:sakov08} (after a spin-up period
of 200 time steps and 50 assimilation cycles). All simulations are started from the same initial ensemble
and use identical sets of observations.

Since the DEnKF has been compared with EnKF and ESRF by \cite{sr:sakov08} and showed the best performance
of all tested methods for this test problem, we only perform
a comparison between the new formulations CEnKF-I (based on (\ref{continuousKalman4b}) with
$\Delta s = 1/4$), CEnKF-II (based on (\ref{continuousKalman4c}) with $\Delta s = 1/4$) and DEnKF.
The mean RMS errors of all three methods can be found in  Tables \ref{table1}  to \ref{table3}.
In Figure \ref{fig_qg1}, we display the RMS error for an optimally chosen 
inflation factor $\delta$ as a function of the localization radius $r_0$. Curves are based on 
the data presented in Tables \ref{table1} to \ref{table3}. We conclude from Figure \ref{fig_qg1} 
that all three filters display a nearly identical performance for an optimal 
$(\delta,r_0)$ parameter choice and that DEnKF shows a slightly better performance for the largest localization radius 
$r_0 = 35$. A similar observation was made by \cite{sr:sakov08} with regard to a comparison between DEnKF 
and ESRF. The differt results for $r_0 = 35$ could be due to the built-in overestimation of the analyzed ensemble
covariance matrix \citep{sr:sakov08}. 

It should be noted that 
CEnKF-II is the least computational expensive of the three methods considered in this study. 


\section{Conclusions and further extensions}

Schur-product-based localization of covariance matrices has become a popular and powerful tool to make
ensemble Kalman filters perform well even under small ensemble sizes. In this note, we have proposed
a new approach to implement Schur-product-based localization seamlessly within the framework of ensemble
Kalman filters. Our approach is based on the formulation of the Kalman update step as 
differential equations in terms of its ensemble members. We have demonstrated for the Lorenz-96 model 
that the resulting methods outperform EnKF with perturbed observations and perform as well as 
standard implementations of ensemble Kalman filters such as ESRF with serial processing  of observations 
and the recently proposed DEnKF. We also implemented a QG model and found that our methods perform
nearly as well as DEnKF which is currently the best available method for this model problem. 
From a computational point of view, the formulation (\ref{continuousKalman4c})
is particularly appealing since it leads to very efficient implementations without the need of
matrix inversions (except when the error covariance matrix ${\bf R}$ is not diagonal) and only
a single evaluation of the ensemble generated covariance matrix $\widetilde{\bf H P}$.

We now outline an number of possible extensions of the formulation (\ref{gradient}).

First we note that (\ref{continuousKalman3})  
can be used in connection with any cost functional $S({\bf x})$ and, hence, 
provides a straightforward generalization of EnKF to nonlinear observation 
operators ${\bf y} = {\bf h}({\bf x})$, i.e.,
\begin{equation} \label{cost2}
S({\bf x}) = 
\frac{1}{2} \left({\bf h}({\bf x}) -{\bf y}\right)^T
{\bf R}^{-1} \left({\bf h}({\bf x}) -{\bf y}\right) .
\end{equation}

Second, as for other localization techniques, the formulation (\ref{continuousKalman4a}) leads to updates in the
ensemble deviations ${\bf X}_a'$ which lie outside the space $\mathbb{S}_f$ in general and, hence, may 
introduce imbalance into the analyzed ensemble members ${\bf x}_i$. 
It seems feasible to restore balance within the proposed framework 
by introducing additional cost functions $S_{\rm pseudo}({\bf x})$
into the formulations (\ref{continuousKalman4a}) or (\ref{continuousKalman4b}), respectively. 
For example, we might require that the divergence of a velocity field ${\bf v}$ remains
``small'' by including a cost functional
\begin{equation} \label{pseudoexample1}
S_{\rm pseudo} = \frac{1}{2r}
\int_\Omega \left(\nabla \cdot {\bf v}\right)^2 \, {\rm d} V,
\end{equation}
where $r>0$ is an appropriate constant. Hence we would modify (\ref{continuousKalman3}) to
\begin{equation} \label{pseudoKalman1}
\frac{\rm d}{{\rm d}s} {\bf x}_i = 
- \frac{1}{2} {\bf P} \left[ \nabla_{{\bf x}_i} \left\{ S({\bf x}_i) + 
S_{\rm pseudo}({\bf x}_i) \right\} 
+ \nabla_{\overline{\bf x}} \left\{ S(\overline{\bf x}) + S_{\rm pseudo}(\overline{\bf x})
\right\} \right] .
\end{equation}

Third, we have focused on deterministic
ensemble Kalman filter formulations in this paper. However, EnKF with perturbed observations \citep{sr:burgers98} can
also be put into the framework of continuous updates and leads naturally to the formulation
\begin{equation}
\frac{\rm d}{{\rm d}s} {\bf x}_i = 
- {\bf P} {\bf H}^T {\bf R}^{-1} \left\{  {\bf H}{\bf x}_i - {\bf y}_i  \right\},
\end{equation}
where ${\bf y}_i$ are now stochastically perturbed observations \citep{sr:burgers98}. Alternatively, we may 
consider the stochastic differential equation
\begin{equation}
{\rm d} {\bf x}_i = 
- {\bf P} {\bf H}^T {\bf R}^{-1} \left\{ {\bf H} {\bf x}_i{\rm d}s - {\bf y}{\rm d}s + 
{\bf R}^{1/2}{\rm d}{\bf w}_i \right\}
\end{equation}
in the ensemble members, where ${\bf w}_i(s) \in \mathbb{R}^k$ denotes standard $k$-dimensional 
Brownian motion. See, for example, \cite{sr:gardiner} for an introduction to stochastic differential equations.

Fourth, we have extensively discussed the continuous formulation of a single Kalman filter analysis step 
for a set of observations given at some time instance $t_j$. We now come back to
the complete ensemble Kalman filter formulation for sequences of
observations at time instances $t_j$, $j=1,\ldots,M$, and intermediate propagation of
the ensemble under the dynamics (\ref{ode}). The continuous formulation
of the ensemble Kalman filter step allows for the following concise
formulation in terms of a differential equation 
\begin{equation} \label{odeEnKF}
\dot{\bf x}_i =
f({\bf x}_i) - \sum_{j=1}^M
\delta(t-t_j)\,{\bf P} \nabla_{{\bf x}_i} {\cal V}_j({\bf X})
\end{equation}
in each ensemble member, where $\delta(\cdot)$ denotes the standard Dirac delta function and
${\cal V}_j({\bf X})$ is the potential (\ref{potential}) with
$S({\bf x})$ replaced by
\begin{equation} 
S_j({\bf x}) = \frac{1}{2} \left({\bf H}{\bf x} -{\bf y}(t_j)\right)^T
{\bf R}^{-1} \left({\bf H}{\bf x} -{\bf y}(t_j)\right).
\end{equation}
One may view (\ref{odeEnKF}) as the original ODE (\ref{ode}) 
driven by a sequence of impluse like contributions due to observations. 
Numerically, it makes sense to regularize these impulses and to 
replace (\ref{odeEnKF}) by 
\begin{equation} \label{mollyEnKF}
\dot{\bf x}_i =
f({\bf x}_i) -  \sum_{j=1}^M
\delta_\varepsilon(t-t_j)\,{\bf P} \nabla_{{\bf x}_i} {\cal V}_j({\bf X}),
\end{equation}
where 
\[
\delta_\varepsilon(s) = \frac{1}{\varepsilon} \psi(s/\varepsilon),
\]
$\psi(s)$ is the standard hat function (or a mollifier in the sense of \cite{sr:friedrichs}),
and $\varepsilon>0$ is an appropriate parameter. Formulation (\ref{mollyEnKF}) can be solved numerically by any
standard ODE solver. There is no longer a strict 
separation between ensemble propagation and filtering. Of course,
the ODE (\ref{mollyEnKF}) becomes extremely stiff as $\varepsilon
\to 0$. A sensible choice is $\varepsilon \sim \Delta t$, where
$\Delta t$ is the natural step-size for the ODE (\ref{ode}).
The mollified formulation (\ref{mollyEnKF}) might be of particular interest in the context of 
the assimilation of non-synoptic measurements, e.g.,
measurements which are arriving semi-continuously in time (see \cite{sr:evensen} for
a discussion of alternative approaches).


\bibliographystyle{plainnat}
\bibliography{survey}


\end{document}